\documentclass[12pt]{article}
\usepackage{amsmath}
\usepackage{amssymb}
\usepackage{amsthm}
\usepackage{amscd}
\usepackage{amsxtra}
\usepackage{verbatim}
\usepackage{xcolor}
\usepackage{color}
\usepackage{enumerate}
\usepackage{mathrsfs}
\usepackage[all]{xy} 

\usepackage{array,arydshln}
\usepackage{bm}

\allowdisplaybreaks

\numberwithin{equation}{section}

\definecolor{dblue}{rgb}{0,0,0.45}
\definecolor{red}{rgb}{0.7,0,0}

 \RequirePackage{geometry}
\newtheorem{theorem}{Theorem}[section]

\newtheorem*{lemma*}{Lemma}

\newtheorem{proposition}[theorem]{Proposition}

\theoremstyle{definition}

\newtheorem{remark}[theorem]{Remark}
\newtheorem{definition}[theorem]{Definition}

\theoremstyle{remark}

\newcommand{\N}{{\mathbb N}}
\newcommand{\R}{{\mathbb R}}

\newcommand{\cC}{{\mathcal C}}

\newcommand{\nn}{\nonumber}
\newcommand{\ve}{\varepsilon}

\newcommand{\vertiii}[1]{{\left\vert\kern-0.25ex\left\vert\kern-0.25ex\left\vert #1 
    \right\vert\kern-0.25ex\right\vert\kern-0.25ex\right\vert}}

\date{}
\begin{document}

\title{Moderate deviations for rough differential equations 
}
\author{  
 Yuzuru \textsc{Inahama}, \quad
  Yong \textsc{Xu}, \quad Xiaoyu \textsc{Yang}
  }
\maketitle

\begin{abstract} 
Small noise problems are quite important for all types of 
stochastic differential equations.
In this paper we focus on rough differential equations
driven by scaled fractional Brownian rough path with 
Hurst parameter $H\in (1/4, 1/2]$. 
We prove a moderate deviation principle for this equation as
the scale parameter tends to zero. 
\end{abstract}

\section{Introduction}

Consider the following stochastic differential equation (SDE)
with a deterministic initial point $a \in \R^e$
driven by a $d$-dimensional standard Brownian motion $(w_t)_{t\in [0,1]}$
scaled by a small parameter  $\ve \in (0,1]$:
\[
dY^{\ve}_t = b(Y^{\ve}_t) dt + \ve\sigma(Y^{\ve}_t) \star dw_t, 
\qquad 
Y^{\ve}_0 = a.
\]
Here, the coefficients $\sigma\colon  \R^e\to  \R^{e\times d}$ and 
$b\colon  \R^e\to  \R^{e}$ are sufficiently regular functions
and 
$\star dw_t$ denotes either the It\^o 
stochastic differential $dw_t$ or the Stratonovich one $\circ dw_t$.

Investigating various limiting behaviors of $Y^\ve =(Y^\ve_t)_{t\in [0,1]}$
as $\ve\searrow 0$ is quite important 
not just for the standard SDE as above but also for many variants of SDEs.
These problems are called small noise problems.
One of the most typical examples are Freidlin-Wentzell's large 
deviation principle (LDP) for $\{Y^\ve\}_{\ve \in (0,1]}$. 
Another example could be a central limit-type theorem 
for $(Y^\ve - Y^0)/\ve$, which states that this process converges in law
to a Gaussian process.

In this paper we take up a moderate deviation principle (MDP),
which is in fact an LDP for $\{Z^\ve\}_{\ve \in (0,1]}$ by definition, 
where we set
\[
Z^{\ve}_t = \frac{Y^{\ve}_t- Y^{0}_t}{\ve^\lambda}, \qquad 0<\lambda <1.
\]
This is equivalent to Freidlin-Wentzell's LDP when $\lambda =0$,
while $\{Z^\ve\}$ satisfies the central limit-type theorem when $\lambda =1$.
Therefore, the MDP bridge the gap between these two famous
limit theorems.

The following is a partial list of
preceding works on MDPs of this kind.
MDPs for various stochastic systems such as
jump-type SDEs \cite{bdg, bdbook},
SDEs with delay  \cite{stz},
stochastic Hamiltonian systems \cite{xgr}, 
slow-fast systems  \cite{gui, gss, ms, hll},
and Volterra-type SDEs \cite{lwyz, jp} have already been proved.
For MDPs for stochastic PDEs, see \cite{wz, wzz, lwz} among others.
In these works, the driving noises are standard, i.e. either Brownian 
or Poisson type.

Study of MDPs for SDEs driven by a (mixed) fractional 
Brownian motion is still in its infancy.
To our knowledge, there are only three works \cite{bds, fyy, yix}.
All of them are quite recent and study the case where Hurst parameter is
larger than $1/2$.
MDPs of this type is not known in the setting of rough path theory.
(Before finishing this work, however, the author was informed of \cite{guz},
in which an MDP is proved for certain rough partial differential equations.
These equations look quite different from those in this paper.)

However, to the best of the authors' knowledge, 
no such result is known for rough differential equations (RDEs)
of standard type.
Our main result (Theorem \ref{thm.LDP})
is an MDP for RDEs driven by a scaled 
fractional Brownian rough path with Hurst parameter $H \in (1/4, 1/2]$.
To prove it, we only use Lyons' continuity theorem,
a Schilder-type LDP for fractional Brownian rough path
and
the contraction principle for LDPs.

The rest of this paper is structured as follows. 
In Section 2, we discuss RDE for the process $Z^\ve$.
Everything in this section is deterministic.
The drift term of the RDE is unbounded, but thanks to \cite{rs}, 
we can make sure that solutions never explode.
Once non-explosion is confirmed, we can show that 
$Z^\ve$ satisfies Lyons' continuity theorem, that is,
it depends continuously on both the driving rough path 
and the small parameter $\ve$. 
Section 3 is a probabilistic part. 
We start by recalling a Schilder-type LDP for fractional 
Brownian rough path on the geometric rough path space.
Our main result is Theorem \ref{thm.LDP}, in which the MDP 
 is stated and proved. 
 The proof is almost immediate from the continuity theorem for $Z^\ve$ 
 since we can combine the contraction principle and the Schilder-type LDP.
 Besides, a central limit-type theorem is also provided in Proposition \ref{prop.CLT}.

%
\medskip

\noindent
{\bf Notation:}~In this paper 
we will use the following notation (unless otherwise specified).
We write $\N =\{1,2, \ldots\}$.
The time interval of (rough) paths 
and stochastic processes is $[0,1]$.
All the vector spaces are over $\R$.

Now we will introduce the notation for some Banach spaces.
(Below, $d, e \in \N$ and $\nabla$ is the standard gradient on 
a Euclidean space.)
\begin{itemize} 
\item
For brevity, we write 
$\R^{e \times d}$
 for the set of real $e \times d$-matrices. 
 The identity matrix of size $e$ is denoted by $\mathrm{Id}_e$ or 
 simply $\mathrm{Id}$.
 Similarly, we write $\R^{e + d}$ for $\R^d\oplus \R^e$.

\item
The set of all continuous path $\varphi\colon [0,1] \to\R^d$ 
is denoted by $\cC (\R^d)$. 
Equipped with the usual sup-norm $\|\varphi\|_{\infty}$,
this is a Banach space.
For $\alpha \in (0,1]$, the set of $\alpha$-H\"older continuous paths
is denoted by $\cC^\alpha (\R^d) := \{ \varphi \in \cC (\R^d)\colon \|\varphi\|_\alpha <\infty\}$, where $\|\varphi\|_\alpha$ is the usual $\alpha$-H\"older seminorm.
Similarly, for $p \in [1,\infty)$, 
the set of  continuous paths of finite $p$-variation is denoted by 
$\cC^{p\textrm{-var}} (\R^d) =\{ \varphi \in \cC (\R^d)
\colon  \|  \varphi\|_{p\textrm{-var}} <\infty\}$,
where $ \|  \varphi\|_{p\textrm{-var}}$ is the usual $p$-variation
 seminorm.
The set of continuous paths that start at $0$ is denoted by $\cC_0 (\R^d)$.
In a similar way, $\cC_0^\alpha (\R^d)$ and $\cC_0^{p\textrm{-var}}$
are defined.

\item
Let $U \subset \R^d$ be a domain.
For $k \in \N \cup \{0\}$,  $C^k (U, \R^e)$ denotes the set of 
$C^k$-functions from $U$ to $\R^e$.
(When $k=0$, we simply write $C (U, \R^e)$ 
instead of $C^0 (U, \R^e)$.)
The set of bounded $C^k$-functions $f \colon U\to \R^e$
whose derivatives up to order $k$ are all bounded 
is denoted by $C_{\mathrm{b}}^k (U, \R^e)$. This is a Banach space with 
the norm $\| f\|_{C_{\mathrm{b}}^k } := \sum_{i=0}^k \|\nabla^i f\|_{\infty}$.
(Here, $ \|\cdot\|_{\infty}$ stands for the usual sup-norm on $U$.)
As usual, we set $C^\infty (U, \R^e) := \cap_{k =0}^\infty C^k (U, \R^e)$
and $C^\infty_{\mathrm{b}} (U, \R^e) := \cap_{k=0}^\infty C^k_{\mathrm{b}} (U, \R^e)$.

\item
Let  $U \subset \R^d$ be a domain
and $\gamma >0$. We write $\gamma = k+\alpha$ 
for $k \in \N$ and $\alpha \in (0,1]$ in a unique way. 
We say $f \colon U\to \R^e$ is of ${\rm Lip}^\gamma$ if 
$f\in C^k_{\mathrm{b}}(U, \R^e)$ and $\nabla^k f$ is 
$\alpha$-H\"older continuous on $U$. 
The set of all such ${\rm Lip}^\gamma$-functions is 
denoted by ${\rm Lip}^\gamma(U, \R^e)$.
The ${\rm Lip}^\gamma$-norm is defined by
\[
\|f\|_{{\rm Lip}^\gamma} := \| f\|_{C_{\mathrm{b}}^k }
+ \sup_{x,y\in U, x\neq y} \frac{|f(x)- f(y)|}{|x-y|^\alpha}.
\]
Note that for $C^k_{\mathrm{b}}(U, \R^e)\subsetneq
{\rm Lip}^k (U, \R^e)$ for every $k\in\N$.

\item
Let $\alpha =1/p\in (1/2,1]$ and $N\in\N$. 
If $w$ belongs to $\cC_0^\alpha (\R^d)$ or $\cC_0^{p\textrm{-var}} (\R^d)$, then we can define 
\[
S_N(w)^m_{s,t} := \int_{0\le t_1 \le \cdots\le t_m\le 1}  dw_{t_1} \otimes\cdots\otimes dw_{t_m},
\qquad 
0 \le s \le t \le 1
\]
as an iterated Young integral for all $m~(1\le m\le N)$. 
We call $S_N(w)$ the natural lift of $w$.

\item
Let $\alpha \in (1/4, 1/2]$. 
We denote by $G\Omega_{\alpha} (\R^d)$ the $\alpha$-H\"older 
geometric rough path space over $\R^d$.
(See \cite{fvbook, lclbook} for a precise definition.)
By definition, $G\Omega_{\alpha} (\R^d)$ is the closure 
of $\{S_{\lfloor 1/\alpha\rfloor}(w) \colon w \in  \cC_0^1 (\R^d)\}$
with respect to the $\alpha$-H\"older rough path metric.
It also coincides with the closure 
of $\{S_{\lfloor 1/\alpha\rfloor}(w) \colon w \in  \cC_0^\beta (\R^d)\}$
for every $\beta \in (1/2,1]$.

 \end{itemize}

%

\section{Deterministic Part}

Let $\alpha \in (1/4, 1/2]$ and $\ve\in (0,1]$. In this section, we consider the following 
rough differential equation (RDE) driven by $\mathbf{x}\in 
G\Omega_{\alpha} (\R^d)$:
\begin{align}  \label{def.RDE_det}
dy^{\ve}_t = b(y^{\ve}_t) dt + \ve\sigma(y^{\ve}_t) d \mathbf{x}_t, 
\qquad 
y^{\ve}_0 = a \in \R^e.
\end{align}
In this work $a$ is arbitrary, but basically fixed. 
In sprit $\ve$ is a small constant. We will let $\ve$ tend to $0$ later. 
It should be recalled that a unique solution of an RDE continously 
depends on both the driving rough path and 
the coefficients with respect to appropriate topologies
under natural assumptions.

Though there are several formulations of RDEs,
we adopt one in Friz-Victoir' book \cite{fvbook} in this paper
because two main preceding results we use  
are both proved in that formulation 
(\cite[Theorem 12.10]{fvbook} and \cite[Theorem 3.1]{rs}). 
In this formulation, a solution of an RDE is a
continuous path in the usual sense and has no ``higher level" objects.
(In any formulation, the first level path of a solution, 
i.e. the component that plays the role of a usual path,  coincides with
 a solution in the above sense after an adjustment of the initial value.)


Let $\kappa \colon (0,1]\to (0,\infty)$ be a
continuous, non-increasing function such that 
$\lim_{\ve\searrow 0} \ve \kappa (\ve)=0$.
(In what follows we understand $0\kappa (0)=0$.)
We are interested in the following object:
\begin{equation} \label{def.Z}
z^{\ve}_t := \frac{y^{\ve}_t- y^{0}_t}{\ve \kappa (\ve)}
\end{equation}
At least formally, one can easily check that
$z^\ve$ satisfies  
\[
d z^{\ve}_t 
  = \left(
               \int_0^1     
                 \nabla b (y^0_t + \theta \ve \kappa (\ve) z^{\ve}_t )  
                 \langle  z^{\ve}_t\rangle d\theta \right) dt 
                 + \kappa( \ve)^{-1}
                   \sigma(y^{0}_t+ \ve \kappa (\ve) z^{\ve}_t) d \mathbf{x}_t, 
\quad z^{\ve}_0 = 0.
\]
The above heuristic consideration leads us to study the following system of RDEs:
\begin{align}
dy^{0}_t &= b(y^{0}_t) dt, 
&y^{\ve}_0 &= a, 
\label{def.RDE_sys1}
\\
d \hat{z}^{\ve}_t 
  &= \left(
               \int_0^1     
                 \nabla b (y^0_t + \theta \ve \kappa (\ve) \hat{z}^{\ve}_t )  
                 \langle  \hat{z}^{\ve}_t\rangle d\theta \right) dt 
                 + 
                   \sigma(y^{0}_t+ \ve \kappa (\ve) \hat{z}^{\ve}_t) d \mathbf{x}_t, 
&\hat{z}^{\ve}_0 &= 0.
\label{def.RDE_sys2}
\end{align}
For the rest of this section, we will show some deterministic properties of 
\eqref{def.RDE_sys1}--\eqref{def.RDE_sys2}.
Note that this system of RDEs makes sense even when $\ve =0$.

\begin{proposition} \label{prop.1}
Let $\alpha \in (1/4, 1/2]$, $\ve\in [0,1]$ and consider the system \eqref{def.RDE_sys1}--\eqref{def.RDE_sys2}  of RDEs driven by $\mathbf{x}\in G\Omega_{\alpha} (\R^d)$.
\\
\noindent 
{\rm (i)}~
Suppose that $\sigma$ is of ${\rm Lip}^{\gamma +1}$ for some
 $\gamma > \alpha^{-1}$ and $b$ is of ${\rm Lip}^2$.
Then, \eqref{def.RDE_sys1}--\eqref{def.RDE_sys2} has  a unique (time-global) solution 
$(y^0, \hat{z}^{\ve})$ for every $\mathbf{x} \in G\Omega_{\alpha} (\R^d)$,
$a\in \R^e$ and $\ve \in [0,1]$.
Moreover, for every $r>0$ and $a\in \R^e$, 
there exists   a constant $C_{a, r}>0$ such that 
 \begin{equation}\label{ineq.RS}
 \|y^0\|_\infty + \|\hat{z}^\ve\|_\infty \le C_{a, r}
  \end{equation}              
  for every $\ve \in [0,1]$ and $\mathbf{x}$ with 
   $\sum_{i=1}^{\lfloor 1/\alpha \rfloor} 
 \|\mathbf{x}^i \|_{\alpha}^{1/i} \le r$.
 Here, $C_{a, r}$ depends only on $r$ and $|a|$
 (and $\sigma$, $b$, $\alpha$).
\\
{\rm (ii)}~Suppose that $\sigma$ is of ${\rm Lip}^{\gamma +1}$ for some
 $\gamma > \alpha^{-1}$ and $b$ is of $C^{3}_{{\rm b}}$.
  Then, 
   \begin{equation}\label{eq.Phi.conti}
    [0,1] \times G\Omega_{\alpha} (\R^d) 
\ni \,\,
(\ve, \mathbf{x})  \mapsto \hat{z}^{\ve}  \,\, \in \mathcal{C}^{\alpha} (\R^e)
\end{equation}
is continuous. 
\end{proposition}

\begin{proof} 
We set 
 $B_R=\{(x, y)\in \R^{e+e} : |x|^2 + |y|^2<R^2 \}$ for $R>0$.
We write $\gamma = m+\lambda$ for a unique 
$(m, \lambda) \in \N\times (0,1]$.
Recall that $\ve \mapsto \ve \kappa(\ve)$ is continuous on $[0,1]$.

First, we show {\rm (i)}.
The diffusion coefficient of the system of RDEs is
\begin{equation}\label{def.sigtilde}
\R^{e+e}\ni 
\begin{pmatrix}
y\\z
\end{pmatrix} 
  \mapsto   
\begin{pmatrix}
{\bf 0}\\
\sigma(y+ \ve \kappa (\ve) z)
\end{pmatrix} 
   \in \R^{(e+e)\times d}
\end{equation}
which is again of ${\rm Lip}^{\gamma +1}$.
Since $\ve \kappa (\ve)$ is bounded in $\ve$, 
its ${\rm Lip}^{\gamma +1}$-norm is bounded 
by $c \| \sigma\|_{{\rm Lip}^{\gamma+1}}$,
where $c>0$ is a constant independent of $\ve$.
The drift of the system of RDEs is 
\begin{equation}\label{def.btilde}
\R^{e+e}\ni \begin{pmatrix}
y\\z
\end{pmatrix} 
  \mapsto 
  \begin{pmatrix}
    b(y)\\
       \int_0^1     
                 \nabla b (y + \theta \ve \kappa (\ve) z)  
                 \langle  z\rangle d\theta
                 \end{pmatrix}                    \in \R^{e+e},
                 \end{equation}
which is clearly locally Lipschitz continuous.
Moreover, it is of linear growth uniformly in $\ve \in [0,1]$, that is, 
\[
|b(y)|+ \left| \int_0^1     
                 \nabla b (y + \theta \ve \kappa (\ve) z)  
                 \langle  z\rangle d\theta \right| 
                 \le 
                  \|b\|_\infty +  \|\nabla b\|_\infty |z|,
                   \qquad y, z\in \R^e.
                                    \]                                  
Now, we use \cite[Theorem 3.1]{rs}, in which Lyons' continuity theorem 
was extended to the case of RDEs with drift vector field of linear growth.   
It assures the existence of a unique global solution 
$(y^0, \hat{z}^{\ve})$ for every $\mathbf{x}$ and $\ve$.    
Inequality \eqref{ineq.RS} is also proved in \cite{rs}.

Next, we show {\rm (ii)}. We write
$\tilde\sigma_\ve (y, z):=\sigma(y+ \ve \kappa (\ve) z)$.
Since it holds for all $y, z\in \R^e$ and  $i~(0 \le i \le m)$ that
\begin{align}  
\lefteqn{
\nabla^i \sigma(y+ \ve \kappa (\ve) z)
- 
\nabla^i \sigma(y+ \ve_0 \kappa (\ve_0) z)
}
\nn\\
&=
\{ \ve \kappa (\ve) - \ve_0 \kappa (\ve_0)\} 
 \int_0^1 d\tau
   \nabla^{i+1} \sigma \left(\tau (y+ \ve \kappa (\ve) z) 
   + (1-\tau) (y+ \ve_0 \kappa (\ve_0) z) \right) 
     \langle  z \rangle,
       \nn
\end{align}
we can easily see that
\[
\lim_{\ve \to \ve_0} \sup_{(y,z)\in B_R} 
 | \nabla^i \sigma(y+ \ve \kappa (\ve) z)
- 
\nabla^i \sigma(y+ \ve_0 \kappa (\ve_0) z)| =0
\]
for all $R>0$ and $i~(0 \le i \le m)$.
Moreover, since 
\begin{align}  
\lefteqn{
\nabla^m \sigma(y_1+ \ve \kappa (\ve) z_1)
- 
\nabla^m \sigma(y_2+ \ve \kappa (\ve) z_2)
}
\nn\\
&=
\int_0^1 d\tau
     \nabla^{m+1} \sigma \left(\tau (y_1+ \ve \kappa (\ve) z_1) 
   + (1-\tau) (y_2+ \ve \kappa (\ve) z_2) \right) 
     \langle  (y_1-y_2) +  \ve \kappa (\ve) (z_1-z_2)\rangle,
     \nn
     \end{align}
we can easily show for all $R>0$ that 
the $\lambda$-H\"older norm on $B_R$ of  
\[
B_R \ni \,\, (y,z) \mapsto \nabla^m \sigma(y+ \ve \kappa (\ve) z) 
-
\nabla^m \sigma(y+ \ve_0 \kappa (\ve_0) z)
\]
converges to $0$ as $\ve \to \ve_0$.
It should be noted that 
we used above only the dominated convergence theorem 
and that $\nabla^{m+1} \sigma$ is bounded and 
uniformly continuouson every bounded subset.
(In other words, 
the H\"older continuity of $\nabla^{m+1} \sigma$ was not used).
Combining these, we can see that 
\[
[0,1]\ni \,\, \ve  \mapsto 
  \tilde\sigma_\ve \,\, \in {\rm Lip}^{\gamma} (B_R, \R^{e\times d})
    \]
is continuous for all $R>0$.

Similarly,  we set $\tilde\beta_\ve (y, z):=\int_0^1     
                 \nabla b (y + \theta \ve \kappa (\ve) z)  
                 \langle  z\rangle d\theta$. 
Then,  essentially in the same way as above, 
we can also show that
\begin{equation}\label{eq.0705-1}
[0,1]\ni \,\, \ve  \mapsto 
  \tilde\beta_\ve \,\, \in {\rm Lip}^{1+\delta} (B_R, \R^e)
  \end{equation}
for every $R>0$ and sufficiently small $\delta >0$.

Now, we use 
\cite[Theorem 12.10 and Remark 12.7 {\rm (i)}]{fvbook},
which is a version of Lyons' continuity theorem for RDEs with drift.
It claims that a solution of such an RDE continuously depends on 
both the driving rough path and the coefficients. 
Thanks to {\rm (i)}, we can use a standard cut-off technique. 
Combining these, we can show that
\begin{equation}
 [0,1]\times
 \{\mathbf{x} \in G\Omega_{\alpha} (\R^d) : \sum_{i=1}^{\lfloor 1/\alpha \rfloor} \|\mathbf{x}^i \|_{\alpha}^{1/i} \le r\}
\ni \,\,
(\ve, \mathbf{x})  \mapsto (y^0, \hat{z}^{\ve})  \,\, \in \mathcal{C}^{\alpha} (\R^{e+e})
\label{eq.0321}
\end{equation}
is continuous for all $r>0$. This proves {\rm (ii)}.

Before closing, we elaborate on the above cut-off procedure.
Set $R =C_{a,r}$ for every $r>0$, where $C_{a,r}>0$ 
is as in \eqref{ineq.RS}.
Choose a smooth function $\chi_r  \colon \R^{e+e} \to [0, 1]$
such that $\chi_r\equiv 1$ on $B_{R}$ and 
$\chi_r\equiv 0$ outside $B_{R+1}$ .
We cut off the coefficients of \eqref{def.RDE_sys2} by 
replacing $\tilde\sigma_\ve (y, z)$ and 
$\tilde\beta_\ve (y, z)$ by $\chi_r (y,z)\tilde\sigma_\ve (y, z)$ and 
$\chi_r (y,z)\tilde\beta_\ve (y, z)$, respectively.
 (No need to cut off the coefficient $b$ in  \eqref{def.RDE_sys1}.)
Suppose that $\sum_{i=1}^{\lfloor 1/\alpha \rfloor} \|\mathbf{x}^i \|_{\alpha}^{1/i} \le r$. Then, due to {\rm (i)},
$(y^0, \hat{z}^{\ve})$ also solves a new system of
RDEs with the new coefficients $\chi_r \tilde\beta_\ve$, $\chi_r \tilde\sigma_\ve$ and $b$.
Since 
$\ve  \mapsto 
 \chi_r \tilde\sigma_\ve  \in {\rm Lip}^{\gamma} (\R^{e+e}, \R^{e\times d})$ and  $\ve  \mapsto 
  \tilde\beta_\ve  \in {\rm Lip}^{1+\delta} (\R^{e+e}, \R^e)
$
 are both continuous, we can apply the above-mentioned facts in 
 \cite{fvbook} to the new system of RDEs to prove 
 the continuity of the map in \eqref{eq.0321}.
\end{proof}

\begin{definition} \label{def.Phi}\eqref{eq.Phi.conti}
We denote by $\Phi\colon [0,1] \times G\Omega_{\alpha} (\R^d) 
\to \mathcal{C}^{\alpha} (\R^e)$ the map defined by \eqref{eq.Phi.conti},
 namely, $\Phi (\ve, \mathbf{x}) =\hat{z}^{\ve}$.
\end{definition}

\begin{proposition} \label{prop.2}
Let $\alpha \in (1/4, 1/2]$ and 
$\ve \in (0,1]$. Let $y^{\ve}$ be a unique solution of 
 RDE \eqref{def.RDE_det} and set $z^{\ve}$ by \eqref{def.Z}.
 Then, we have 
 \[
 z^{\ve} = \Phi (\ve, \kappa(\ve)^{-1} \mathbf{x}), 
 \qquad  
    \mathbf{x}\in G\Omega_{\alpha} (\R^d), \, \ve \in (0,1].
     \]
      Here, $\kappa(\ve)^{-1} \mathbf{x}$ is the dilation of $ \mathbf{x}$ by $\kappa(\ve)^{-1}>0$.
\end{proposition}

\begin{proof} 
Let $x \in \mathcal{C}^{1} (\R^{d})$ and denote its natural lift by
$\mathbf{x}:= S_{\lfloor 1/\alpha \rfloor} (x)$. 
In this case, $y^{\ve}$ is a unique solution of 
the following Riemann-Stieltjes ODE:
\[
dy^{\ve}_t = b(y^{\ve}_t) dt + \ve\sigma(y^{\ve}_t) dx_t, 
\qquad 
y^{\ve}_0 = a \in \R^e.
\]
We can see from this that
\begin{align}  
z^{\ve}_t 
 &=  \ve^{-1} \kappa (\ve)^{-1} \left\{
    \int_0^t \{ b(y^{\ve}_s) - b(y^{0}_s) \}ds  
        + \ve
                        \int_0^t \sigma(y^{\ve}_s) dx_s 
                          \right\}
                          \nn\\
                           &=
                            \ve^{-1} \kappa (\ve)^{-1}
    \int_0^t \{ b(y^{0}_s +   \ve \kappa (\ve)z^{\ve}_s ) -b(y^{0}_s) \}ds  
        + \kappa (\ve)^{-1}
                        \int_0^t \sigma(y^{0}_s +   \ve \kappa (\ve)z^{\ve}_s) dx_s
                          \nn\\
                          &=
                           \int_0^t \left(
                            \int_0^1 
                             \nabla b(y^{0}_s+   \theta\ve \kappa (\ve)z^{\ve}_s)  \langle z^{\ve}_s \rangle d\theta            \right)                                             ds
                                   + \kappa (\ve)^{-1}
                        \int_0^t \sigma(y^{0}_s +   \ve \kappa (\ve)z^{\ve}_s) dx_s.
                       \nn                       
                          \end{align}
Hence, we have
$z^{\ve} = \Phi (\ve, \kappa (\ve)^{-1} S_{\lfloor 1/\alpha \rfloor} (x))$ in this case.
For a general $\mathbf{x} \in G\Omega_{\alpha} (\R^d)$,
we take $\{ x_k\}_{k\in \N} \subset \mathcal{C}^{1} (\R^{d})$ such 
that $\lim_{k\to\infty}S_{\lfloor 1/\alpha \rfloor} (x_k) =\mathbf{x}$ 
in $G\Omega_{\alpha} (\R^d)$
and use the continuity of $\Phi (\ve, \cdot)$ and $\mathbf{x} \mapsto y^{\ve}$
for each fixed $\ve \in (0,1]$.
\end{proof}

%

\section{Probabilistic Part}
In this section, we take  parameters as follows.
Let $H \in (1/4, 1/2]$.
If $H \in (1/3, 1/2]$, we take $\alpha\in (1/3, H)$.
if $H \in (1/4, 1/3]$, we take $\alpha\in (1/4, H)$.
Note that $\lfloor H^{-1} \rfloor=\lfloor\alpha^{-1} \rfloor$.)

Denote by $(w^H_t)_{t\in [0,1]} =(w^{H,1}_t, \ldots, w^{H,d}_t)_{t\in [0,1]}$
be a $d$-dimensional fractional Brownian motion with Hurst parameter $H$.
A canonical rough path lift of $w^H$ is denoted by
$\mathbf{W}^H$ and is called  fractional Brownian rough path 
with Hurst parameter $H$.
It is viewed as a $G\Omega_{\alpha} (\R^d)$-valued random variable.
For $m \in\N$, we denote by $w^H(m)$ be a piecewise linear approximation of $w^H$
associated with $\{i/2^m : 0\le i \le 2^m\}$. 
It is known that $S_{\lfloor 1/\alpha \rfloor} (w^H(m))$ converges 
(at least) in probability to $\mathbf{W}^H$ with respect to the $\alpha$-H\"older
rough path topology.

We denote by $\mathcal{H}^{H}(\R^{d})$ be the Cameron-Martin space of $w^H$.
Each $h\in \mathcal{H}^{H}(\R^{d})$ is $H$-H\"older continuous
and of finite $\{H +(1/2)\}^{-1}$-variation (see \cite{fvjfa, fggr}).
Note that $1\le \{H +(1/2)\}^{-1}< 4/3$.  
Hence, $S_{\lfloor 1/\alpha \rfloor} (h)$ is well-defined in the variation setting,
although its H\"older regularity is not so clear a priori. 
However, it is known that
$S_{\lfloor 1/\alpha \rfloor} (h) \in G\Omega_{\alpha} (\R^d)$ (see \cite{fvjfa}).
The injection 
$S_{\lfloor 1/\alpha \rfloor}\colon \mathcal{H}^{H} (\R^{d})\hookrightarrow G\Omega_{\alpha} (\R^d)$ is locally Lipschitz continuous.

Let $\ve\mathbf{W}^H$ be the dilation of $\mathbf{W}^H$ by $\ve \in (0,1]$.
A Schilder-type LDP is known, that is, 
$\{\ve\mathbf{W}^H\}_{\ve \in (0,1]}$ satisfies 
an LDP on $G\Omega_{\alpha} (\R^d)$ 
as $\ve\searrow 0$ with speed $\ve^{-2}$ and  a good rate function $J$,
which is defined by 
\[
J (\mathbf{x})  = \left\{
\begin{array}{ll}
\|h\|_{\mathcal{H}^{H}(\R^{d})}^2 /2
   & (\mbox{if $\mathbf{x} =S_{\lfloor 1/\alpha \rfloor} (h)$ for 
     some $h\in \mathcal{H}^{H}(\R^{d})$}), \\
 +\infty & (\mbox{otherwise}).
\end{array}
\right.
\]
(See \cite[Theorem 15.55]{fvbook}.)
Moreover, $\{\ve\mathbf{W}^H\}_{\ve \in (0,1]}$ is exponentially tight 
on $G\Omega_{\alpha} (\R^d)$,
due to a Fernique-type theorem for $\mathbf{W}^H$. 
(See \cite[Theorem 15.33]{fvbook}.).

For $\ve \in [0,1]$, let $Y^\ve$ be a unique solution of \eqref{def.RDE_det} with 
$\mathbf{x}$ being replaced by $\mathbf{W}^H$, namely, 
\begin{align}  \label{def.RDE_random}
dY^{\ve}_t = b(Y^{\ve}_t) dt + \ve\sigma(Y^{\ve}_t) d \mathbf{W}^H_t, 
\qquad 
Y^{\ve}_0 = a \in \R^e,
\end{align}
and set for $\ve \in (0,1]$
\begin{equation} \label{def.randomZ}
Z^{\ve}_t := \frac{Y^{\ve}_t- Y^{0}_t}{\ve \kappa (\ve)}.
\end{equation}
Clearly, $Z^{\ve}$ is a $\mathcal{C}^{\alpha} (\R^{d})$-valued random variable.

\begin{proposition}\label{prop.CLT}
Consider the case $\kappa \equiv 1$.
Suppose that $\sigma$ is of ${\rm Lip}^{\gamma +1}$ for some
 $\gamma > H^{-1}$ and $b$ is of $C^{3}_{{\rm b}}$.
Then,  as $\ve\searrow 0$, 
\[
Z^{\ve}_t = \frac{Y^{\ve}_t- Y^{0}_t}{\ve} \,\, \to\,\,  \Phi (0, \mathbf{W}^H)
\]
in $\mathcal{C}^{\alpha} (\R^{d})$ almost surely.
Moreover, $\Phi (0, \mathbf{W}^H)$ is a mean-zero Gaussian process.
\end{proposition}

\begin{proof} 
Let $H^{-1}< \alpha^{-1} < \gamma \wedge (\lfloor H^{-1}\rfloor +1)$.
The convergence is immediate from Propositions \ref{prop.1} and \ref{prop.2}
(and Definition \ref{def.Phi}, too). 
So, it remains to show that Gaussian property. 
Since 
\[
\Phi (0, \mathbf{W}^H) = \lim_{m\to\infty} \Phi (0, S_{\lfloor 1/\alpha \rfloor} (w^H(m)))
\qquad \mbox{a.s.,}  
\]
it suffices to check that $\Phi (0, S_{\lfloor 1/\alpha \rfloor} (w^H(m)))$, 
which will be denoted by $\Xi (m)$, is Gaussian with mean zero.  
By definition, $\Xi (m)$ solves the following Riemann-Stieltjes ODE:
\[
d \Xi (m)_t 
  = 
                 \nabla b (y^0_t)  
                 \langle   \Xi (m)_t\rangle  dt 
                 + 
                   \sigma(y^{0}_t) d w^H(m)_t, 
\qquad \Xi (m)_0 = 0.
\]
Let $M$ be a unique solution of the following $e\times e$ matrix-valued ODE:
\[
dM_t = \nabla b (y^0_t)  M_t dt, \qquad  M_0 = {\rm Id}_e.
\]
Note that $\nabla b$ is viewed as an $e\times e$ matrix-valued 
function. Then, $M_t$ is invertible and non-random and we have
\[
\Xi (m)_t 
  = M_t \int_0^t M_s^{-1}\sigma(y^{0}_s)  w^H(m)^{\prime}_s ds
      \]
for all $t\in [0,1]$.
Note that, for all $s$, $w^H(m)^{\prime}_s$ can be written as a linear combination 
of $\{w^H_{i/2^m} \colon 0\le i \le 2^m\}$. 
So, $\Xi (m)_t$ can be written as a limit of  linear combinations
of $\{w^H_{s} \colon 0\le s \le 1\}$, which implies that $\Xi (m)$ is a 
mean-zero Gaussian process.
\end{proof}

Now we provide our main theorem.
It is an MDP for RDEs driven by 
fractional Brownian rough path with Hurst parameter $H\in (1/4, 1/2]$.
A prominent example of $\kappa$ is $\kappa (\ve)=\ve^{-\theta}$ 
for $0<\theta <1$.

\begin{theorem}\label{thm.LDP}
Let $H \in (1/4, 1/2]$ and $\alpha \in (0,H)$.
Suppose that $\kappa \colon (0,1]\to (0,\infty)$ is a
continuous, non-increasing function such that 
$\lim_{\ve\searrow 0} \kappa (\ve)= +\infty$ and
$\lim_{\ve\searrow 0} \ve \kappa (\ve)=0$.
Suppose further that $\sigma$ is of ${\rm Lip}^{\gamma +1}$ for some
 $\gamma > H^{-1}$ and $b$ is of $C^{3}_{{\rm b}}$.

  Then, $\{Z^{\ve}\}_{\ve\in (0,1]}$ satisfies an LDP
  in $\mathcal{C}^{\alpha} (\R^{d})$  as $\ve\searrow 0$ 
  with speed $\kappa (\ve)^2$
  and a good rate function $I$ given by
  \[
I (\xi) = \inf \{  \|h\|_{\mathcal{H}^{H}(\R^{d})}^2 /2 \colon
  \mbox{$h\in \mathcal{H}^{H}(\R^{d})$ such that 
  $\xi =  \Xi^h$}\}, 
  \quad \xi \in \mathcal{C}^{\alpha} (\R^{d}).
  \]
  As usual we set $\inf \emptyset =+\infty$. 
    Here, $\Xi^h$ stands for a unique solution of the following Young ODE
    driven by $h$: 
    \begin{equation}\label{eq.0705-2}
    d\Xi^h_t = 
                 \nabla b (y^0_t)  
                 \langle   \Xi^h_t\rangle  dt 
                 + 
                   \sigma(y^{0}_t) dh_t, 
\qquad \Xi^h_0 = 0.
        \end{equation}
          \end{theorem}

\begin{proof} 
The larger $\alpha$ is, the stronger the claim of the theorem becomes.
Hence, it is enough to assume  
$H^{-1}< \alpha^{-1} < \gamma \wedge (\lfloor H^{-1}\rfloor +1)$.

Consider the family of point masses $\{\delta_\ve\}_{\ve\in (0,1]}$ on $[0,1]$.
Clearly, it satisfies an LDP on  $[0,1]$
as $\ve\searrow 0$ with speed $\kappa (\ve)^2$ and  a good rate function $K$,
where $K (0) := 0$ and $K(s) :=+\infty$ if $0<s\le 1$.
It is also clear that 
$\{\delta_\ve\}_{\ve \in (0,1]}$ is exponentially tight on $[0,1]$.

By a general fact for LDPs for product measures
(see \cite[p. 129]{dzbook} for instance), 
 $\{(\ve, \kappa(\ve)^{-1}\mathbf{W}^H)\}_{\ve\in (0,1]}$  satisfies 
an LDP on $[0,1]\times G\Omega_{\alpha} (\R^d)$ 
as $\ve\searrow 0$ with speed $\kappa (\ve)^2$ and  a good rate function $\hat{J}$,
where 
\[
\hat{J} (\ve, \mathbf{x})  := \left\{
\begin{array}{ll}
\|h\|_{\mathcal{H}^{H}(\R^{d})}^2 /2
   & (\mbox{if $\ve =0$ and $\mathbf{x} =S_{\lfloor 1/\alpha \rfloor} (h)$ for 
     some $h\in \mathcal{H}^{H}(\R^{d})$}), \\
 +\infty & (\mbox{otherwise}).
\end{array}
\right.
\]

By Proposition \ref{prop.2}, we have
$Z^{\ve} = \Phi (\ve, \kappa(\ve)^{-1}\mathbf{W}^H)$
and $\Phi$ is continuous by Proposition \ref{prop.1}.
Therefore, we can use the contraction principle \cite[Theorem 4.2.1]{dzbook} to 
obtain the desired LDP for $\{Z^\ve\}_{\ve\in (0,1]}$ with
a good rate function $I$  given as follows:
\[
I (\xi) = \inf \{  \|h\|_{\mathcal{H}^{H}(\R^{d})}^2 /2 \colon
  \mbox{$h\in \mathcal{H}^{H}(\R^{d})$ such that 
  $\xi =   \Phi (0,  S_{\lfloor 1/\alpha \rfloor} (h))$}\}, 
  \quad \xi \in \mathcal{C}^{\alpha} (\R^{d}).
  \]
Noting that $\Xi^h = \Phi (0,  S_{\lfloor 1/\alpha \rfloor} (h))$,  we completes the proof.
 \end{proof}

\begin{remark} \label{rem.BMcase}
By specializing $H=1/2$ in Theorem \ref{thm.LDP}, we recover 
known moderate deviation results for usual SDEs at least to some extent.
(Since we use rough path theory, the conditions on $b$ and $\sigma$   
in this remark are stronger than those in preceding works.)
Note that $W^H$ is Stratonovich-type 
Brownian rough path in this case.
In this remark, $\kappa$ is the same as in Theorem \ref{thm.LDP}.

(1)~Suppose that $\sigma$ is of ${\rm Lip}^{\gamma +1}$ for some $\gamma > 2$ and $b$ is of $C^{3}_{{\rm b}}$.
Then, the solution $Y^\ve$ of RDE \eqref{def.RDE_random}
coincides with a unique solution of 
 the following usual Strotonovich-type SDE driven by 
standard Brownian motion $(w^{1/2}_t)_{t\in [0,1]}$:
\[
dy^{\ve}_t = b(y^{\ve}_t) dt + \ve\sigma(y^{\ve}_t) \circ dw^{1/2}_t, 
\qquad 
y^{\ve}_0 = a \in \R^e.
\]
Note that this SDE has a unique (non-exploding) solution
because both  $b$ and $\sigma$ are of $C^{2}_{{\rm b}}$
(see a corollary in \cite[p. 106]{iwbook} for instance).
Therefore, an MDP  for $\{y^{\ve} \}_{\ve \in (0,1]}$
(i.e.  an LDP for $\{z^{\ve} \}_{\ve \in (0,1]}$,
where $z^{\ve} := (y^{\ve}- y^{0})/ \{\ve \kappa (\ve)\}$) 
is 
a special case of Theorem \ref{thm.LDP} above. 

(2)~Next, we discuss It\^o-type SDEs.
Suppose that $\sigma$ is of $C^{4}_{{\rm b}}$ and $b$ is of $C^{3}_{{\rm b}}$.
Instead of \eqref{def.RDE_random}, we consider the following RDE: 
\begin{align}  \label{def.RDE_ito}
d\tilde{Y}^{\ve}_t = \tilde{b}_{\ve}(\tilde{Y}^{\ve}_t) dt + \ve\sigma(\tilde{Y}^{\ve}_t) d \mathbf{W}^{1/2}_t, 
\qquad 
Y^{\ve}_0 = a \in \R^e,
\end{align}
where 
\begin{equation}\nn
\tilde{b}_{\ve}^i (y) :=b^i (y) - \frac{\ve^2}{2}\sum_{j=1}^d 
\sum_{k=1}^e \sigma_{kj}(y)  \cdot
 \partial_k \sigma_{ij}(y) ,
\qquad
 y \in\R^e, \,\, 1\le i \le e.
\end{equation}
In other words, $b$ in  \eqref{def.RDE_random}
was replaced by $\tilde{b}_{\ve}$, which is of $C^{3}_{{\rm b}}$ again.
Note $\ve^2$ in front of the It\^o-Stratonovich correction term.
Then, $\tilde{Y}^\ve$ coincides with a unique solution of 
 the following usual It\^o-type SDE:
 \[
d\tilde{y}^{\ve}_t = b(\tilde{y}^{\ve}_t) dt + \ve\sigma(\tilde{y}^{\ve}_t)  dw^{1/2}_t, 
\qquad 
\tilde{y}^{\ve}_0 = a \in \R^e.
\]
It should be noted that since $\tilde{b}_{\ve}$ depends on $\ve$,
an MDP for $\{\tilde{Y}^{\ve} \}_{\ve \in (0,1]}$
(equivalently, those for $\{\tilde{y}^{\ve} \}_{\ve \in (0,1]}$)
is not proved in Theorem \ref{thm.LDP}.
But, we can slightly modify Theorem \ref{thm.LDP} to cover this case
as follows.
First, $b$ in RDE \eqref{def.RDE_sys2} is replaced by $b _\ve$.
Then,  the (new) drift vector field of this RDE, as 
a function of $\ve$, 
still satisfies the same property as in \eqref{eq.0705-1}.
Moreover, 
thanks to the factor $\ve^2$,
the limiting skeleton ODE \eqref{eq.0705-2} remains unchanged, i.e.
the correction term vanishes from this ODE.
For these reasons, we can see the same MDP holds for $\{\tilde{Y}^{\ve} \}_{\ve \in (0,1]}$ (and $\{\tilde{y}^{\ve} \}_{\ve \in (0,1]}$), too.
\end{remark}

\medskip

\noindent
{\bf Acknowledgments}:~
This work was partly supported by JSPS KAKENHI (Grant No. 20H01807), 
the Key International (Regional)
Cooperative Research Projects of the NSF of China (Grant 12120101002) and the NSF of China (Grant 12072264).
%
%

%

\bigskip
\begin{flushleft}
  \begin{tabular}{ll}
    Yuzuru \textsc{Inahama}
    \\
    Faculty of Mathematics,
    \\
    Kyushu University,
    \\
    744 Motooka, Nishi-ku, Fukuoka, 819-0395, JAPAN.
    \\
    Email: {\tt inahama@math.kyushu-u.ac.jp}
  \end{tabular}
\end{flushleft}

\begin{flushleft}
  \begin{tabular}{ll}
    Yong \textsc{Xu}
    \\
School of Mathematics and Statistics,
    \\
Northwestern Polytechnical University, 
    \\
Xi'an, 710072, CHINA.
    \\
    Email: {\tt hsux3@nwpu.edu.cn}
  \end{tabular}
\end{flushleft}

\begin{flushleft}
  \begin{tabular}{ll}
    Xiaoyu \textsc{Yang}
    \\
    Faculty of Mathematics,
    \\
    Kyushu University,
    \\
    744 Motooka, Nishi-ku, Fukuoka, 819-0395, JAPAN.
    \\
    Email: {\tt yangxiaoyu@yahoo.com}
  \end{tabular}
\end{flushleft}


\begin{thebibliography}{00}

\bibitem{bds}
Bourguin, S.; Dang, T.; Spiliopoulos, K.;
Moderate deviation principle for multiscale systems driven by fractional Brownian motion.
To appear in J. Theor. Probab. (2023). 
https://doi.org/10.1007/s10959-023-01235-y

\bibitem{bdbook}
Budhiraja, A.; Dupuis, P.;
{\it Analysis and approximation of rare events.}  Springer, 2019.


\bibitem{bdg}
Budhiraja, A.; Dupuis, P.; Ganguly, A.;
Moderate deviation principles for stochastic differential equations with jumps.
Ann. Probab. 44 (2016), no.3, 1723--1775.

\bibitem{dzbook}
 Dembo, A.,  Zeitouni, O.;
{\it Large deviations techniques and applications.} Springer, 2009.


\bibitem{fggr}
Friz, P.; Gess, B.; Gulisashvili, A.; Riedel, S.;
The Jain-Monrad criterion for rough paths and applications to random Fourier series and non-Markovian H\"ormander theory. 
Ann. Probab. 44 (2016), no. 1, 684--738. 



\bibitem{fvjfa}
Friz, P. K.; Victoir, N. B.; 
A variation embedding theorem and applications.
J. Funct. Anal. 239 (2006), no.2, 631--637.

\bibitem{fvbook}
Friz, P. K.; Victoir, N. B.;
 {\it Multidimensional stochastic processes as rough paths.} Cambridge University Press, Cambridge, 2010.


\bibitem{fyy}
Fan, X.; Yu, T.; Yuan, C.;
Asymptotic behaviors for distribution dependent SDEs driven by fractional Brownian motions.
Preprint (2022). arXiv: 2207.01525.

\bibitem{gss}
Gesteratos, I.; Salins, M.; Spiliopoulos, K.;
Moderate deviations for systems of slow-fast stochastic reaction-diffusion equations.
Preprint (2021). arXiv: 2101.00085.


\bibitem{gui}
Guillin, A.;
Averaging principle of SDE with small diffusion: moderate deviations.
Ann. Probab.31(2003), no.1, 413--443.




\bibitem{guz}
Gussetti, E.;
Pathwise central limit theorem and moderate deviations via rough paths for SPDEs with multiplicative noise. 
Preprint (2023). arXiv: 2307.10965.  


\bibitem{hll}
Hong, W.; Li, G.; Li, S.;
Multi-Scale McKean-Vlasov SDEs: Moderate Deviation Principle in Different Regimes.
Preprint (2023). arXiv:2306.11569.

\bibitem{iwbook}
Ikeda, N.; Watanabe, S.;
 {\it Stochastic differential equations and diffusion processes. }
Second edition. 
North-Holland Publishing Co., Amsterdam; Kodansha, Ltd., Tokyo, 1989.



\bibitem{jp}
Jacquier, A.; Pannier, A.;
Large and moderate deviations for stochastic Volterra systems.
Stochastic Process. Appl.149 (2022), 142--187.

\bibitem{js}
Jacquier, A.; Spiliopoulos, K.;
Pathwise moderate deviations for option pricing.
Math. Finance 30 (2020), no. 2, 426--463.

\bibitem{ms}
Morse, M. R.; Spiliopoulos, K.;
Moderate deviations for systems of slow-fast diffusions.
Asymptot. Anal. 105 (2017), no.3-4, 97--135.




\bibitem{lwyz}
Li, Y.; Wang, R.; Yao, N.; Zhang, S.;
A moderate deviation principle for stochastic Volterra equation.
Statist. Probab. Lett. 122 (2017), 79--85.

\bibitem{lwz}
Li, Y.; Wang, R.;  Zhang, S.;
Moderate deviations for a stochastic heat equation with spatially correlated noise.
Acta Appl. Math. 139 (2015), 59--80.

\bibitem{lclbook}
Lyons, T.; Caruana, M.; L\'evy, T.;
{\it Differential equations driven by rough paths. }
  Lecture Notes in Math., 1908. Springer, Berlin, 2007.


\bibitem{rs}
Riedel, S.; Scheutzow, M.;
Rough differential equations with unbounded drift term.
J. Differential Equations 262 (2017), no.1, 283--312.


\bibitem{stz}
Suo, Y.; Tao, J.; Zhang, W.; 
Moderate deviation and central limit theorem for stochastic differential delay equations with polynomial growth.
Front. Math. China13(2018), no.4, 913--933.


\bibitem{wz}
Wang, R.; Zhang, T.;
Moderate deviations for stochastic reaction-diffusion equations with multiplicative noise.
Potential Anal. 42 (2015), no.1, 99--113.

\bibitem{wzz}
Wang, R.; Zhai, J.; Zhang, T.;
A moderate deviation principle for 2-D stochastic Navier-Stokes equations.
J. Differential Equations 258 (2015), no.10, 3363--3390.


\bibitem{xgr}
Xu, J.; Gong, J.; Ren, J.;
A moderate deviation principle for stochastic Hamiltonian systems.
ESAIM Probab. Stat. 27 (2023), 558--575.

\bibitem{yix}
Yang, X.; Inahama, Y.; Xu, Y.;
Moderate deviations for two-time scale system with mixed fractional
Brownian motion.
In preparation (2023). 

\end{thebibliography}
\end{document}